\documentclass[11pt]{amsart}
\usepackage{amsmath,amsthm, amscd, amssymb, amsfonts}
\usepackage[all]{xy}
\usepackage{mathrsfs}
\usepackage{enumitem}
\usepackage{ wasysym}
\usepackage{hyperref}

\usepackage[ansinew]{inputenc}
\usepackage{graphicx,fancyhdr}

\usepackage[dvips, dvipsnames, usenames]{color}

\newcommand{\comment}[1]{}

\numberwithin{equation}{subsection}
\newtheorem{theorem}{Theorem}
\newtheorem{lemma}[theorem]{Lemma}
\newtheorem{coro}[theorem]{Corollary}

\newtheorem{prop}[theorem]{Proposition}

\theoremstyle{definition}
\newtheorem{definition}[theorem]{Definition}
\newtheorem{example}[theorem]{Example}

\theoremstyle{remark}
\newtheorem{remark}[theorem]{Remark}

\def\pf{\begin{proof}}
\def\epf{\end{proof}}

\newcommand{\ku}{ \Bbbk}

\newcommand\I{\mathbb I}

\newcommand\N{\mathbb N}

\newcommand\Z{\mathbb Z}

\newcommand\D{\mathcal{D}}

\newcommand{\cP}{\mathcal{P}}

\newcommand{\Ss}{{\mathcal S}}
\newcommand\T{\mathcal{T}}
\newcommand\cV{\mathcal{V}}

\newcommand{\Ug}{\mathfrak U}
\newcommand{\ujor}{{\Ug}^{\text{\texttt{jordan}}}}
\newcommand{\wujor}{{\widetilde\Ug}^{\mathtt{jordan}}}

\newcommand\Der{\operatorname{Der}}
\newcommand\End{\operatorname{End}}
\newcommand\id{\operatorname{id}}
\newcommand\gr{\operatorname{gr}}
\newcommand\GK{\operatorname{GKdim}}
\newcommand{\Hom}{\operatorname{Hom}}

\newcommand{\yd}[1]{{}^{ #1 }_{ #1 }\mathcal{YD}}

\newcommand{\toba}{\mathscr{B}}
\newcommand{\ot}{\otimes}

\newcommand{\planck}{\saturn}

\newcounter{tabla}\stepcounter{tabla}

\begin{document}

\title[Corrigendum: Liftings of the Jordan plane]{Corrigendum: Liftings of the Jordan plane}

\author[Andruskiewitsch; Angiono; Heckenberger]
{Nicol\'as Andruskiewitsch, Iv\'an Angiono, Istv\'an Heckenberger}

\address{FaMAF-CIEM (CONICET), Universidad Nacional de C\'ordoba,
Medina A\-llen\-de s/n, Ciudad Universitaria, 5000 C\' ordoba, Rep\'
ublica Argentina.} \email{nicolas.andruskiewitsch@unc.edu.ar, ivan.angiono@unc.edu.ar}

\address{Philipps-Universität Marburg,
Fachbereich Mathematik und Informatik,\newline
Hans-Meerwein-Straße,
D-35032 Marburg, Germany.} \email{heckenberger@mathematik.uni-marburg.de}

\thanks{
\newline  N. A. and I. A. were partially supported by CONICET,
FONCyT-ANPCyT, Secyt (UNC). The work of N. A. was partially done during a visit to the University of Clermont. }

\begin{abstract}
We complete the classification of the pointed Hopf algebras with finite Gelfand-Kirillov dimension that are 
liftings of  the  Jordan plane over a nilpotent-by-finite group, correcting the statement in \cite{AAH}.
\end{abstract}

\maketitle

\setcounter{tocdepth}{1}
\renewcommand{\thesubsection}{\arabic{subsection}}

 \subsection*{Introduction}
In the paper \cite{AAH} we stated the classification of the pointed Hopf algebras with finite Gelfand-Kirillov dimension that are 
liftings of  either the  Jordan plane or the super Jordan plane over a nilpotent-by-finite group.
But we overlooked one possibility, namely to deform degree one relations and therefore the 
classification in \emph{loc. cit.} of liftings of Jordan planes is not complete. Here we fill the gap.
It turns out that the missed example is essentially a Hopf algebra introduced by C. Ohn in 1992, see \cite{O}.

\medbreak  Throughout $\ku$ is an algebraically closed field of characteristic $0$.
Recall that   $\cV(1, 2)$ is the braided vector space with basis $x_1$, $x_2$ and braiding $c$
given by
$c(x_i \otimes x_1) =  x_1 \otimes x_i$, $c(x_i \otimes x_2) = ( x_1 + x_2) \otimes x_i$, $i  = 1,2$.
Here is the revised version of \cite[Proposition 4.2]{AAH}.

\begin{prop}\label{prop:lifting-jordan}
Let $G$ be a nilpotent-by-finite group and  let $H$ be a pointed Hopf algebra with finite $\GK$  such that 
\begin{itemize}[leftmargin=*]\renewcommand{\labelitemi}{$\circ$}
\item $G(H) \simeq G$ and

\item the infinitesimal braiding of $H$ is isomorphic to $\cV(1, 2)$.
\end{itemize}
Then there exists a Jordanian YD-triple $\D=(g, \chi, \eta)$ for $\Bbbk G$ such that  either

\begin{enumerate}[leftmargin=*,label=\rm{(\Roman*)}]

\item\label{item:ss} $H\simeq \Ug(\D,0)$ or $H\simeq \Ug(\D,1)$, introduced in \cite[\S 4.1]{AAH}; or

\item\label{item:ss-indesc} $\chi=\varepsilon$ and there exists $\xi\in \Der_{\varepsilon,\varepsilon}(\ku G, \ku)$, $\xi\ne 0$,
such that $H \simeq \Ug_{\xi}(\D,0)$ or $H\simeq \Ug_{\xi}(\D,1)$ see Definition \ref{def:A}; or 

\item\label{item:indesc} $\chi=\varepsilon$ and  $H\simeq \ujor(\D)$, see Definition \ref{def:B}.
\end{enumerate}
\end{prop}

Conversely, any of these Hopf algebras is pointed and has finite $\GK$, actually  $\GK \Bbbk G +2$.
See Lemmas \ref{lema:A}, \ref{lema:B} and  \cite[Proposition 4.2]{AAH}.
Notice that if $\chi=\varepsilon$ and $\xi = 0$, then   $\Ug_{0}(\D,\lambda) \simeq \Ug(\D, \lambda)$, introduced in \cite[\S 4.1]{AAH}.

\medbreak
The subspace of $(g,1)$ skew-primitive elements in a  Hopf algebra in case \ref{item:ss} is decomposable as $G$-module,
while in \ref{item:ss-indesc} is decomposable as $\langle g \rangle$-module but it is an indecomposable $G$-module,
and in \ref{item:indesc} it is an indecomposable $\langle g \rangle$-module. Thus Hopf algebras from different cases could not be isomorphic. Whether Hopf algebras in the same case are isomorphic is treated as in \cite[\S 4.1]{AAH}.

\medbreak
This note is organized as follows. 
In Section \ref{subsec:Ujord}  the minimal Hopf algebra missing in \cite[Proposition 4.2]{AAH} and its relation with \cite{O} are described.
In Section \ref{subsec:gap} we discuss  the gap.
 Proposition \ref{prop:lifting-jordan} is proved in Section \ref{subsec:proof}.

\subsection*{Notation} We keep the notations from \cite{AAH}.
Let $G$ be a group, let $\ku G$ be its group algebra and let $\widehat G$ be its  group of characters.
Given  $\chi \in\widehat G $, recall that
\begin{align*}
\Der_{\chi,\chi}(\ku G, \ku) &= \{\eta\in (\ku G)^*: \eta(ht) = \chi(h)\eta(t) + \chi(t)\eta(h) \quad\forall h,t \in G\}.
\end{align*}

A  collection $\D=(g, \chi, \eta)\in Z(G) \times \widehat G \times \Der_{\chi,\chi}(\ku G, \ku)$ is a 
\emph{YD-triple} for $\ku G$ if $\eta(g) = 1$.
Then the vector space $\cV_g(\chi,\eta)$ with a basis $(x_i)_{i\in\I_2}$ belongs to $\yd{\ku G}$,
with the coaction  $\delta(x_i) = g\otimes x_i$, $i\in \I_2$, and the action given by
\begin{align*}
h\cdot x_1 &= \chi(h) x_1,& h\cdot x_2&=\chi(h) x_2 + \eta(h)x_{1},&  h&\in \ku G.
\end{align*}
When $\chi(g) = 1$ we say that  $\D=(g, \chi, \eta)$ is a \emph{Jordanian} YD-triple.

\medbreak
 Let $L$ be a Hopf algebra. The $\Delta$,  $\varepsilon$ and $\Ss$ denote respectively the comultiplication,
 the counit and the antipode.
 The  group of group-like elements is denoted by $G(L)$. Also
 the space of $(g,h$)-primitive elements is $\cP_{g,h}(L) =\{\ell \in L: \Delta(\ell) = \ell\otimes h + g \otimes \ell\}$, 
 where $g,h \in G(L)$,   and $\cP(L) = \cP_{1,1}(L)$ is the space of primitive elements.
The adjoint action of $G(L)$ on $L$ is denoted by $g\cdot \ell := g\ell g^{-1}$, $g\in G(L)$, $\ell\in L$. 

\subsection*{Competing interests} The authors declare none.

 \subsection{The Jordanian enveloping algebra of  $s\ell(2)$}\label{subsec:Ujord}
 Let $\wujor$ be the algebra generated by $a_1, a_2, g, g^{-1} $ with defining relations
 \begin{align}\label{eq: U-jordan-rels1}
 g^{\pm1}  g^{\mp1} &= 1, & ga_1 &= a_1g + (g-g^2), &  ga_2 &= a_2g + a_1g.
 \end{align}
 It is easy to see that $\wujor$ is a Hopf algebra by imposing $g \in G(\wujor)$ and $a_1,a_2 \in \cP_{g, 1}(\wujor)$.
 We introduce
 \begin{align}\label{eq:def-z}
 z &= a_1a_2 - a_2a_1 -\dfrac{a_1^2}{2} + a_2 +\dfrac{1}{2}a_1 \in \wujor
 \end{align}
 \begin{lemma}\label{lema:z-skewprim}
The element $z$ belongs to $\cP_{g^2,1}(\wujor)$
and commutes with $g$.
\end{lemma}

\pf We compute
\begin{align*}
\Delta (z) &= a_1a_2\otimes 1 + a_1g \otimes a_2 + g a_2 \otimes a_1 + g^2 \otimes  a_1a_2 \\
& -a_2a_1\otimes 1 - a_2g \otimes a_1 - g a_1 \otimes a_2 -  g^2 \otimes  a_2a_1 \\
& -\dfrac{1}{2} a_1^2\otimes 1 -\dfrac{1}{2} \left( a_1g + g a_1\right) \otimes a_1-\dfrac{1}{2}  g^2 \otimes  a_1^2 \\
& + a_2\otimes 1 + g \otimes a_2  +\dfrac{1}{2} a_1\otimes 1 + \dfrac{1}{2}  g \otimes a_1 \\
&= z \otimes 1 + g^2 \otimes z + \left(a_1g  - g a_1 + g -  g^2 \right)\otimes a_2 \\
& + \left(g a_2 - a_2g  -\dfrac{1}{2} ( a_1g + g a_1) + \dfrac{1}{2}  g -  \dfrac{1}{2}  g^2  \right)\otimes a_1\\
&= z \otimes 1 + g^2 \otimes z;
\end{align*}
here 
$g a_2 - a_2g  -\dfrac{1}{2} ( a_1g + g a_1) + \dfrac{1}{2}  g -  \dfrac{1}{2}  g^2 = 
\dfrac{1}{2} \left(a_1g - g a_1 + (g - g^2)\right) = 0$.

It remains to prove that $\gamma(z)=0$, where
$\gamma \in \End_{\ku}(\wujor)$ is given by
$\gamma (x) = gxg^{-1}-x$, for all $ x\in \wujor$.
Note that 
\begin{align*}
\gamma(xy) &=\gamma(x)(\gamma(y)+y)+x\gamma(y)&\text{for all } x,y &\in \wujor.
\end{align*}
From \eqref{eq: U-jordan-rels1} we have that
\begin{align} \label{eq:gamma}
\gamma(a_1)&=1-g, &  \gamma(a_2)&=a_1.
\end{align}
Therefore,
\begin{align*}
\gamma(z)&=\gamma\Big(a_1a_2+\Big(a_2+\frac 12 a_1\Big)(1-a_1)\Big)\\
&=\gamma(a_1)(\gamma(a_2)+a_2)+a_1\gamma(a_2)\\
&+\gamma\Big(a_2+\frac 12a_1\Big)(\gamma(1-a_1)+1-a_1)
+\Big(a_2+\frac 12a_1\Big)\gamma(1-a_1).\\
\intertext{By using \eqref{eq:gamma} we obtain that}
\gamma(z)&=(1-g)(a_2+a_1)+a_1^2\\
&\quad +\Big(a_1+\frac 12 (1-g)\Big)(g-a_1)
+\Big(a_2+\frac 12a_1\Big)(g-1)\\
&=a_2+a_1-(a_2g+2a_1g+g-g^2)+a_1^2
+\Big(a_2+\frac 12a_1\Big)(g-1)\\
&\quad +a_1g-a_1^2+\frac 12(g-g^2)+\frac 12a_1(g-1)+\frac 12(g-g^2).
\end{align*}
Now it follows easily that $\gamma(z)=0$.
\epf

The Jordanian enveloping algebra of  $s\ell(2)$ is 
\begin{align}\label{eq:def-ujordan}
\ujor := \wujor / \langle z\rangle.
\end{align}
By Lemma \ref{lema:z-skewprim}, $\ujor$ is a Hopf algebra quotient of
$\wujor$.
By abuse of notation the images of $g, a_1, a_2$ in $\ujor$
are denoted by the same symbols.

\begin{remark}\label{rem:ujorlambda}
For each $\lambda \in \ku $ let
\begin{align}\label{eq:def-ujordan2}
\ujor_{\lambda} := \wujor / \langle z-\lambda(1-g^2)\rangle.
\end{align}
Then $\ujor_{\lambda}$ is a Hopf algebra, since $z-\lambda(1-g^2)\in 
\cP_{g^2,1}(\wujor)$.

Let us now fix $\lambda ,\mu \in \ku $. Let $U$ be the algebra
$$U =\ku \langle g,g^{-1},a_1,a_2\rangle /\langle gg^{-1}-1,
g^{-1}g-1\rangle. $$
Then $U$ has a unique Hopf algebra
structure such that $g,g^{-1}\in G(U)$ and $a_1,a_2\in \cP_{g,1}(U)$.
Moreover, there exists a well-defined Hopf algebra map
$$ \varphi_{\lambda,\mu}:U\to \ujor_{\lambda}, \quad
g\mapsto g,\, a_1\mapsto a_1,\, a_2\mapsto a_2+\mu(1-g).
$$
It is easily checked that
$$ ga_1-a_1g-g+g^2,\,ga_2-(a_2+a_1)g\in \ker \varphi_{\lambda,\mu}. $$
Moreover, for $z\in U$ defined as in \eqref{eq:def-z} we obtain that
\begin{align*}
\varphi_{\lambda,\mu}(z)-z&=a_1\mu(1-g)-\mu(1-g)a_1+\mu(1-g)=\mu(1-g^2).
\end{align*}
Since $z=\lambda(1-g^2)\in \ujor_{\lambda}$, we conclude that
$\varphi_{\lambda,\mu}$ induces a surjective Hopf algebra map
$$ \varphi_{\lambda,\mu}:\ujor_{\lambda+\mu}\to \ujor_\lambda. $$
It follows that $\varphi_{0,\lambda}:\ujor_{\lambda}\to \ujor $
is a Hopf algebra isomorphism.
\end{remark}

\begin{remark}\label{rem:Ujord}
For any $\planck\in\ku$, the Hopf algebra $U_{\planck}$ was introduced by Christian Ohn in \cite{O};
this is the algebra generated over $\ku$ by $K,Y,T^{\pm 1}$ with relations:
\begin{align}\label{def:Uh1}
TT^{-1} = T^{-1}T = 1, \quad [K,T] &=T^2-1, \quad [Y,T] = -\frac {\planck}{2}(KT+TK),
\\ \label{def:Uh2}
[K,Y] &=-\frac12(YT+TY+YT^{-1}+T^{-1}Y),
\end{align}
with the Hopf algebra structure  of $U_{\planck}$    determined by
 $ T\in G(U_{\planck})$ and $ X, Y \in \cP_{T^{-1}, T}(U_{\planck})$.
It is easy to see that the  the Hopf algebras $U_{\planck}$ with $\planck\neq 0$ are all isomorphic so we fix one of them.
The appellative \emph{Jordanian} was introduced  by Alev and Dumas to the best of our knowledge.
We claim that $\ujor_{\lambda}$ is isomorphic to  the Hopf subalgebra  $\Ug$ of $U_{\planck}$ generated by 
\begin{align}
x &= KT^{-1}, &  y &= YT^{-1}, & g &=T^{-2};
\end{align}
we choose these variables  to have $ x, y \in \cP_{g, 1}(\Ug)$. Now  \eqref{def:Uh1} implies
\begin{align}\label{action} 
g \cdot x &= x+ 2 (1-g), & g \cdot y &= y - 2\planck \left(x + (1- g)\right).
\end{align}
We perform a new change of variables:
\begin{align*}
a_1 &= \frac{1}{2}x, & a_2 &=  -\frac{1}{4\planck}y -  \frac{1}{4}x;
\end{align*}
these new variables satisfy \eqref{eq: U-jordan-rels1}. Now 
\eqref{def:Uh2} translates succesively into 
\begin{align*}
xy-yx & =-2y-{\planck}x^2+\frac {\planck}4(1-g^{2})
\end{align*}
and then into
\begin{align*}
z &= -\frac{1}{32}(1-g^2). 
\end{align*}
 That is, $\Ug \simeq \ujor_{-\frac{1}{32}}$.
\end{remark}

\begin{remark}
The algebra $U_{\planck}$  can be described as an iterated Ore extension:
\begin{equation}U_{\planck}=\ku[T^{\pm}][x\,;\delta][y\,;\,\sigma,D] 
\end{equation}
with $\delta$ a derivation of $\ku[T^{\pm}]$, $\sigma$ an automorphism of $\ku[T^{\pm}][x\,;\delta]$ and 
$D$ a $\sigma$-derivation of $\ku[T^{\pm}][x\,;\delta]$ defined by:
\begin{align}\textstyle xT&=Tx+\underbrace{(T-T^{-1})}_{=\delta(T)}\\
yT&=\underbrace{T}_{=\sigma(T)}y+\underbrace{(-\planck Tx-\frac \planck 2 (T-T^{-1}))}_{=D(T)}\\
yx&=\underbrace{(x+2)}_{=\sigma(x)}y+\underbrace{\planck x^2-\frac \planck 4(1-T^{-4})}_{=D(x)}.
\end{align}
\end{remark}

\begin{prop}\label{prop:Ujord}
There exist a derivation $\delta_1$ of $R:=\Bbbk[g,g^{-1}]$,  a derivation $\delta_2$ of $S:=R[a_1;\id,\delta_1]$
and an automorphism $\sigma$ of $S$
such that $\ujor$ is isomorphic to the Ore extension $S[a_2;\sigma,\delta_2]$.

Hence $\ujor$ is a  noetherian domain of Gelfand-Kirillov 3, and
 the monomials $g^ja_1^{i_1} a_2^{i_2}$ form a PBW-basis  of $\ujor$.
\end{prop}

\pf  We leave the verification of the first claim to the reader as a long but straightforward exercise: the derivations $\delta_1:R\to R$, $\delta_2:S\to S$ satisfy
\begin{align*}
\delta_1(g) &= g^2-g, & \delta_2(g) &= -a_1g, & \delta_2(a_1) &= \frac{1}{2}a_1(1-a_1),
\end{align*}
and $\sigma$ is given by $\sigma(g)=g$, $\sigma(a_1)=a_1+1$. The rest is standard.
\epf

\begin{coro}\label{cor:Ujord}
The  Hopf algebra   $\ujor$ is pointed and $\gr \ujor$ is isomorphic to the bosonization of the Jordan plane
by the group algebra of the infinite cyclic group. \qed
\end{coro}

\subsection{The gap and how to fix it}\label{subsec:gap}
We fix a group $G$.
Let $H$ be a pointed Hopf algebra with coradical filtration $(H_n)_{n\in \N_0}$ such that $ G(H) \simeq G$. 
Then $H_1/ H_0 \simeq V \# \ku G$, where $V \in \yd{\ku G}$  is the infinitesimal braiding of $H$.
For $g\in G$,  the space of $(g,1)$ skew-primitives $\cP_{g,1}(H)$
satisfies 
\[\cP_{g,1}(H) \cap H_0 = \ku (1-g) \text{ and }  \cP_{g,1}(H)  / \left(\cP_{g,1}(H) \cap H_0 \right) \simeq  V_g.\]

\medbreak 
Now assume that  $V \simeq\cV_g(\chi,\eta)$ for a YD-triple $\D = (g, \chi, \eta)$ over $\ku G$.
Thus $V = V_g$ and we have an exact sequence  of  $G$-modules
\begin{align*}
\xymatrix{0 \ar@{->}[r] & \ku (1-g) \ar@{->}[r] &  \cP_{g,1}(H)   \ar@{->}[r] ^{\varpi}&  \cV_g(\chi,\eta) \ar@{->}[r] & 0.}
\end{align*}
Since $g\in Z(G)$, one has  $\ku (1-g) \subset \cP_{g,1}(H)^{\varepsilon}$. Hence  $\chi \neq \varepsilon$ implies that 
\begin{align*}
 \cP_{g,1}(H) \simeq  \ku (1-g)  \oplus \cV_g(\chi,\eta)
\end{align*}
and we have a morphism of Hopf algebras $\pi: \T(\cV_g(\chi,\eta))  \to H$, where
$\T(\cV_g(\chi,\eta)) = T(\cV_g(\chi,\eta)) \# \ku G$.  
In particular the proof of \cite[Prop. 4.3]{AAH} goes over without changes.

We assume for the rest of this Section 
that  the infinitesimal braiding $V$ of $H$ is isomorphic to $\cV_g(\varepsilon,\eta)$ 
for a YD-triple $\D = (g, \varepsilon, \eta)$ as  Yetter-Drinfeld module over $\ku G$.
Under this assumption, $\cP_{g,1}(H)$ might be indecomposable.

\begin{example} The  indecomposability of $\cP_{g,1}(H)$ could happen in other situations. 
Here is a simple example.
Let $A$ be the algebra generated by $a, \gamma^{\pm1}$, where $ \gamma^{-1}$ is the  inverse of $ \gamma$
and the relation  $\gamma a\gamma^{-1} = a + (1- \gamma)$ holds, so that $A$ is not commutative. 
Then $A$ is a pointed Hopf algebra by declaring that $\gamma$ is a group-like and $a$ a
 $(\gamma,1)$ skew-primitive element. Observe that  $\cP_{g,1}(A)$ is indecomposable.
Let $\Gamma \simeq \Z$.
It can be shown that $\gr A \simeq T(V) \otimes \ku \Gamma$, 
where $V$ has dimension 1 and is the infinitesimal braiding of $A$.
But $\cP_{g,1}(A)$ is indecomposable and  there is no surjective morphism of Hopf algebras $T(V) \otimes \ku \Gamma \to A$.
\end{example}

Back to our situation,  let us pick $a_1, a_2 \in  \cP_{g,1}(H)$ such that 
$\varpi(a_j) = x_j$, $j= 1,2$ and set $a_0 = 1-g$.   Then 
there are $\zeta \in  \Der_{\varepsilon,\varepsilon}(\ku G, \ku)$ and a linear map $\xi: \ku G  \to \ku$  such that
the action of $h\in G$ on  $\cP_{g,1}(H)$ is given  in the basis
$(a_0, a_1, a_2)$  by 
\begin{align}\label{eq:action-h-skew-prim}
\rVert h\rVert = \begin{pmatrix}
1 &  \zeta(h) &  \xi(h) \\
0 & 1 & \eta(h) \\
0 &0 & 1 
\end{pmatrix}.
\end{align}
Notice that $\Der_{\varepsilon,\varepsilon}(\ku G, \ku) = \Hom_{\text{gps}}(G, (\ku, +))$
and that $\xi$ is a kind of differential operator of degree 2, meaning that
\begin{align}\label{eq:xi-property}
\xi (hk) &= \xi(h)+\zeta(h)\eta(k) +\xi(k) & &\text{for all }h,k\in G.
\end{align}
Thus if $\zeta\neq 0$, then the claim  \cite[Prop. 4.2, page p. 669, line 8]{AAH} is not true.
To correct  this  we consider the subalgebra $A$  generated by $g$ and $\cP_{g,1}(H)$, 
 a Hopf subalgebra of $H$.
The action of $g$ on  $\cP_{g,1}(H) = \cP_{g,1}(A)$  in the basis
$(a_0, a_1, a_2)$ is given by 
\begin{align}\label{eq:action-g-skew-prim}
\rVert g\rVert  =\begin{pmatrix}
1 &  \zeta(g) &  \xi(g) \\
0 & 1 & 1 \\
0 &0 & 1 
\end{pmatrix}.
\end{align}

As $g\in Z(G)$, we have that $\xi(gh) = \xi(hg)$  for all $h\in G$, so \eqref{eq:xi-property}
says that
\begin{align}\label{eq:relation-zeta-eta}
\zeta(h) &= \eta(h)\zeta(g) & &\text{for all }h\in G.
\end{align}

We consider two cases:

\begin{enumerate}[leftmargin=*,label=\rm{(\Alph*)}]
\item\label{item:A} $\zeta(g) = 0$. 
Then $\zeta=0$ by \eqref{eq:relation-zeta-eta} and $\xi\in \Der_{\varepsilon,\varepsilon}(\ku G, \ku)$ by \eqref{eq:xi-property}.

\item\label{item:B} $t:= \zeta(g) \neq 0$, the Jordanian case. In the basis $(a_0, t^{-1}a_1, t^{-1}a_2 - t^{-2}\xi(g) a_1)$,  the action of $g$    is given by 
$\small \begin{pmatrix}
	1 &  1 &  0 \\
	0 & 1 & 1 \\
	0 &0 & 1 
\end{pmatrix}$.
We still denote the new basis by $(a_0, a_1, a_2)$; that is, we may assume that $\zeta(g)=1$, $\xi(g)=0$. By \eqref{eq:relation-zeta-eta}, $\zeta=\eta$, and by \eqref{eq:xi-property},
$\xi (hk)= \xi(h)+\eta(h)\eta(k) +\xi(k)$ for all $h,k\in G$.  
\end{enumerate}

We shall see that the following Hopf algebras exhaust the  case \ref{item:A}.

\begin{definition} \label{def:A} 
Let $\D = (g, \varepsilon, \eta)$ be a  YD-triple, 
$\xi\in \Der_{\varepsilon,\varepsilon}(\ku G, \ku)$ and $\lambda\in\Bbbk$.
 We define
$\Ug_{\xi}(\D,\lambda)$ as the algebra generated by $h\in G$, $a_1$, $a_2$ 
with defining relations being those of $G$ and
\begin{align}\label{eq:Uxi-def-rel-1}
&ha_1 - a_1 h, & &h\in G;
\\\label{eq:Uxi-def-rel-2}
&ha_2 - \left(a_2+\eta(h)a_1+\xi(h)(1-g)\right)h, & &h\in G;
\\\label{eq:Uxi-def-rel-3}
&a_1a_2 - a_2a_1 -\dfrac{a_1^2}{2}  -\lambda(1-g^2).
\end{align}

\end{definition}

As we said already,    $\Ug_{0}(\D,\lambda)\simeq \Ug(\D,\lambda)$,  introduced in \cite[\S 4.1]{AAH}.

\begin{lemma}\label{lema:A}
$\Ug_{\xi}(\D,\lambda)$  is a Hopf algebra with comultiplication determined by 
\begin{align*}
G(\Ug_{\xi}(\D,\lambda)) &= G& &\text{and}& a_1,a_2 &\in \cP_{g,1}(\Ug_{\xi}(\D,\lambda)).
\end{align*}
 Thus $\Ug_{\xi}(\D,\lambda)$ is pointed. 
 The set $\{a_1^m a_2^n h \, | \, m,n\in \N_0, \, h\in G \}$ is a basis of $\Ug_{\xi}(\D,\lambda)$;
   $\gr \Ug_{\xi}(\D,\lambda)\simeq \toba(\cV(1,2))\# \ku G$ and
 \begin{align*}
 \GK \Ug_{\xi}(\D,\lambda) &= \GK \Bbbk G +2.
 \end{align*}
 In particular, if $G$ is nilpotent-by-finite, then $\GK \Ug_{\xi}(\D,\lambda) < \infty$. 
\end{lemma}

\pf Left to the reader.
\epf

We shall see that the following Hopf algebras exhaust the  case \ref{item:B}.

\begin{definition} \label{def:B} 
Let $\D = (g, \varepsilon, \eta)$ be a  YD-triple and 
define $\xi\in (\ku G)^*$  by  $\xi(h) =\tfrac{1}{2}(\eta(h)^2-\eta(h))$, $h \in G$.
We introduce  $\ujor(\D)$ as the algebra generated by $h\in G$, $a_1$, $a_2$ 
with defining relations those of $G$, \eqref{eq:Uxi-def-rel-2} and
\begin{align}\label{eq:Ujor-def-rel-1}
&ha_1 - \left(a_1 +\eta(h)(1-g) \right) h, & &h\in G.
\\\label{eq:Uxi-def-rel-3bis}
&a_1a_2 - a_2a_1 -\dfrac{a_1^2}{2} + a_2 +\dfrac{1}{2}a_1.
\end{align}
\end{definition}

Observe that $\xi$, needed in \eqref{eq:Uxi-def-rel-2}, satisfies \eqref{eq:xi-property} with $\zeta = \eta$.
The proof of the following Lemma is also standard.

\begin{lemma}\label{lema:B}
$\ujor(\D)$  is a Hopf algebra with structure determined by 
\begin{align*}
G(\ujor(\D)) &= G& &\text{and}& a_1,a_2 &\in \cP_{g,1}(\ujor(\D)).
\end{align*}
Thus $\ujor(\D)$ is pointed.
 The set $\{a_1^m a_2^n h \, | \, m,n\in \N_0, \, h\in G \}$ is a basis of $\ujor(\D)$;
 $\gr \ujor(\D)\simeq \toba(\cV(1,2))\# \ku G$ and
\begin{align*}
\GK \ujor(\D) &= \GK \Bbbk G +2.
\end{align*}
In particular, if $G$ is nilpotent-by-finite, then $\GK \ujor(\D) < \infty$. \qed
\end{lemma}

\subsection{Proof of Proposition \ref{prop:lifting-jordan}} \label{subsec:proof}

Let $G$ be a nilpotent-by-finite group and  let $H$ be a pointed Hopf algebra with finite $\GK$  such that 
$ G(H) \simeq G$ and the infinitesimal braiding $V$ of $H$ is isomorphic to $\cV(1, 2)$.
By \cite[Lemma 2.3]{AAH}, there exists a unique YD-triple $\D = (g, \chi, \eta)$ such that 
$V \simeq\cV_g(\chi,\eta)$ in $\yd{\ku G}$.
By \cite[Lemma 3.7]{AAH}, $\gr H \simeq \toba(\cV(1,2)) \# \ku G$, hence $H$ is generated 
by $\cP_{g,1}(H)$  and $G$ as algebra.

\medbreak
If $\chi \neq \varepsilon$, then the proof of \cite[Prop. 4.1]{AAH}  implies that $H$ is isomorphic
either to $\Ug(\D,0)$ or $\Ug(\D,1)$, the Hopf algebras introduced in \cite[\S 4.1]{AAH}.

\medbreak
Assume that $\chi = \varepsilon$.  Pick a basis $(a_0 = 1-g, a_1, a_2)$ such that
any  $h\in G$  acts on  $\cP_{g,1}(H)$  by
\eqref{eq:action-h-skew-prim} where 
$\zeta \in  \Der_{\varepsilon,\varepsilon}(\ku G, \ku)$ and $\xi\in  (\ku G )^*$ satisfies \eqref{eq:xi-property}.  
Let  $A$ be the subalgebra generated by $\cP_{g,1}(H)$.
As explained above we consider two cases.

\medbreak
\noindent{\bf Case \ref{item:A}:} $\zeta(g) = 0$, thus $\zeta = 0$.
Even if \cite[Proposition 4.2]{AAH} does not apply in general since we may have $\xi\ne 0$, 
it does apply to $A$ up to changing the base to $(a_0, a_1, \widetilde{a}_2)$ where $\widetilde{a}_2 := a_2 - \xi(g) a_1$, see \eqref{eq:action-h-skew-prim}. 
Call the new basis again $(a_0, a_1, a_2)$ by abuse of notation.
Hence $A \simeq \Ug(\D',\lambda)$ where 
$\D' = (g, \chi_{\vert \langle g \rangle}, \eta_{\vert \langle g \rangle})$ is a YD-triple over the subgroup
$\langle g \rangle$ of $G$ and  $\lambda\in \{0,1\}$. In particular the following equality holds in $H$:
$$a_2a_1=a_1a_2-\tfrac{1}{2}a_1^2+\lambda(1-g^2) .$$

\medbreak
We first claim that $A$ is stable under the action of $G$.  
Indeed let $G$ act on the free algebra generated by $g^{\pm 1}$, $a_1$, $a_2$, where $G$ acts trivially on $g$, and by \eqref{eq:action-h-skew-prim} on $a_1$, $a_2$.
As $g$ is central, the action of each $h\in G$ preserves the defining ideal of $A$, so $G$ acts on $A$.

\medbreak
We next claim that $H \simeq A\rtimes\ku G /I$, where $I$ is the ideal that identifies the two copies of $g$
where $\rtimes$ stands for smash product.
Indeed, the inclusions $A\hookrightarrow H$, $\Bbbk G\hookrightarrow H$ induce a Hopf algebra map $\psi:A\rtimes \Bbbk G/ I  \to H$.
As $\gr H\simeq \toba(V)\#\Bbbk G$, $H$ is generated by $a_1$, $a_2$ and $G$, so $\psi$ is surjective.
On the other hand, $(A\rtimes \Bbbk G/ I )_1$ is spanned by the set $\{1\ot h, a_1\ot h, a_2\ot h | h\in G\}$.
The image of this set under $\psi$ is linearly independent, which implies that $\psi_{\vert (A\rtimes \Bbbk G/ I )_1}$ is injective. 
By \cite[5.3.1]{Mo-libro}, $\psi$ is injective, and the claim follows. 
As a consequence, the set $\{a_1^m a_2^n h \, | \, m,n\in \N_0, \, h\in G \}$ is a basis of $H$.

\medbreak
Finally, we see  that there is a Hopf algebra map $\Ug_{\xi}(\D,\lambda) \to H$; since this map sends a basis to a basis, we conclude that $H\simeq \Ug_{\xi}(\D,\lambda)$.

\medbreak
\noindent{\bf Case \ref{item:B}. }$\zeta(g) \neq 0$.
As discussed above, we may assume that $\zeta = \eta$.
Recall that we are assuming that $\GK H < \infty$.
We claim that 

\begin{enumerate}[leftmargin=*,label=\rm{\textbf {(\roman*)}}]
\item\label{item:lifting-subalgebra} There exists a Hopf algebra isomorphism $A \simeq \ujor$, cf. \eqref{eq:def-ujordan}.

\item\label{item:lifting-condition-xi} $\xi(h)=\tfrac{1}{2}(\eta(h)^2-\eta(h))$ for all $h\in G$.

 \item\label{item:lifting-general-case} $A$ is stable under the adjoint action of $G$ and $H \simeq A\rtimes\ku G /I$, where $I$ is the ideal that identifies the two copies of $g$.
 
 \item\label{item:lifting-final-step} The set $\{a_1^m a_2^n h \, | \, m,n\in \N_0, \, h\in G \}$ is a basis of $H$
 and $H\simeq \ujor(\D)$.
\end{enumerate}

\medbreak

\ref{item:lifting-subalgebra}: It is easy to see that
there exists a Hopf algebra surjective map $\widetilde{\pi}:\wujor\to A$, 
which applies $g$, $a_1$, $a_2$ to the corresponding elements of $A$.
Hence $\widetilde{\pi}(z)\in \cP_{g^2,1}(A)$, by Lemma \ref{lema:z-skewprim}.
Now, as $g\ne g^2$ and $\gr H\simeq \toba(V)\# \Bbbk G$, we have that $\cP_{g^2,1}(H)=\cP_{g^2,1}(H) \cap H_0 = \ku (1-g^2)$;
thus there exists $\lambda\in \Bbbk$ such that $\widetilde{\pi}(z)=\lambda(1-g^2)$, which implies that 
$\widetilde{\pi}$ factors through  a map
$\pi: \ujor_{\lambda} \twoheadrightarrow A$. 
The set $\{g^k,a_1g^k, a_2g^k: k\in \Z\}$ is linearly independent in $H$, so $\pi_{\vert (\ujor_{\lambda})_1}$ is injective. By \cite[5.3.1]{Mo-libro}, $\pi$ is an isomorphism. 
Up to composing with $\varphi_{0,\lambda}$, see Remark \ref{rem:ujorlambda}, we may assume that $\lambda=0$.

\medbreak

\ref{item:lifting-condition-xi}: Given $h\in G$, let $\gamma_h\in \End_{\ku} H$ be given by
\begin{align*}
\gamma_h (x)&=h x h^{-1} -x & \text{for all } & x\in H.
\end{align*}
Note that $\gamma_h(xy)=\gamma_h(x)(\gamma_h(y)+y)+x\gamma_h(y)$
for all $x,y\in H$.
From \eqref{eq:action-h-skew-prim},
\begin{align} \label{eq:gamma-on-U}
\gamma_h(a_1)&=\eta(h)(1-g), &  \gamma_h(a_2)&=\eta(h)a_1+\xi(h)(1-g).
\end{align}
Therefore,
\begin{align*}
\gamma_h(z)
&=\eta(h)(1-g)(\eta(h)a_1+\xi(h)(1-g)+a_2)+a_1(\eta(h)a_1+\xi(h)(1-g))\\
&+\Big(\eta(h)a_1+\xi(h)(1-g)+\tfrac 12 \eta(h)(1-g)\Big)(-\eta(h)(1-g)+1-a_1)
\\ & -\Big(a_2+\tfrac 12a_1\Big)\eta(h)(1-g)
=\left(\tfrac 12 \eta(h) -\tfrac 12 \eta(h)^2+\xi(h)\right)\left(1-g^2 \right).
\end{align*}
By \ref{item:lifting-subalgebra}, $z=0$, so $\gamma_h(z)=0$. Thus, $\xi(h)=\tfrac 12\left( \eta(h)^2- \eta(h)\right)$.

\medbreak

\ref{item:lifting-general-case}:  Let  $G$ act on the free algebra generated by $g^{\pm 1}$, $a_1$, $a_2$, where $G$ acts trivially on $g$, and by \eqref{eq:action-h-skew-prim} on $a_1$, $a_2$. Each $h\in G$ fixes the defining relations $gg^{-1} - 1$, $g^{-1}g - 1$, $ga_1 - a_1g - g+g^2$, $z$, and
\begin{align*}
h\cdot & \left( ga_2 - a_2g - a_1g \right) 
= ga_2- a_2g- a_1g+\eta(h) \left(ga_1-a_1g -(1-g)g \right),
\end{align*}
so the action descends to $A$.  The proof of 
\ref{item:lifting-final-step} is as in Case \ref{item:A}. \qed

\subsection*{Acknowledgements} 
N. A. thanks Fran\c cois Dumas for conversations on and help with the computations in this Note.

\end{document}